\documentclass[review]{elsarticle}

\usepackage{multirow}
\usepackage{breakurl}

\usepackage{lineno,hyperref}
\modulolinenumbers[5]

\journal{Journal}

\usepackage{numcompress}

\begin{document}

\begin{frontmatter}

\title{Solution to dynamic economic dispatch with prohibited operating zones via MILP}

\author[mymainaddress,mysecondaryaddress]{Shanshan~Pan}
\address[mymainaddress]{College of Electrical Engineering, Guangxi University,
Nanning 530004, P. R. China}

\author[mysecondaryaddress]{Jinbao~Jian\corref{mycorrespondingauthor}}
\cortext[mycorrespondingauthor]{Corresponding author}
\ead{jianjb@gxu.edu.cn}
\address[mysecondaryaddress]{Guangxi Colleges and Universities Key Laboratory of Complex System Optimization and Big Data Processing, Yulin Normal University, Yulin 537000, P. R. China}

\author[mythirdaryaddress]{Linfeng~Yang}
\address[mythirdaryaddress]{College of Computer Electronics and Information, Guangxi University, Nanning 530004, P. R. China}

\begin{abstract}
Dynamic economic dispatch (DED) problem considering prohibited operating zones (POZ), ramp rate constraints, transmission losses and spinning reserve constraints is a complicated non-linear problem which is difficult to solve efficiently. In this paper, a mixed integer linear programming (MILP) method is proposed to solve such a DED problem. Firstly, a novel MILP formulation for DED problem without considering the transmission losses, denoted by MILP-1, is presented by using perspective cut reformulation technique. When the transmission losses are considered, the quadratic terms in the transmission losses are replaced by their first order Taylor expansions, and then an MILP formulation for DED considering the transmission losses, denoted by MILP-2, is obtained. Based on MILP-1 and MILP-2, an MILP-iteration algorithm (MILP-IA) is proposed to solve the complicated DED problem. The effectiveness of the MILP-1 and MILP-IA are assessed by several cases and the simulation results show that both of them can solve to competitive solutions in a short time.
\end{abstract}


\end{frontmatter}


\section{Introduction}
Dynamic economic dispatch (DED) is a fundamental tool for the optimal economic operation in power system which aims at allocating the customers' load demands among
the available thermal power generating units in an economic, secure and reliable way \cite{2010Xia}. In practice, some thermal or hydro generating units may have prohibited operating zones (POZ) due to the physical limitations of power plant components, e.g., vibrations in a shaft bearing are amplified in a certain operating region \cite{1993Reserve, 1996GA}. The resulting disjoint operating regions lead to a discontinuous generation cost function, which complicates the DED problem.

Comparing to the dynamic economic dispatch with prohibited operating zones (DED-POZ), the static economic dispatch with prohibited operating zones (SED-POZ) which handles only a single load level at a particular time instant is more achievable. Thus, in the past few decades, a myriad of optimization methods have been presented to deal with SED-POZ.
They include genetic algorithm (GA) \cite{1996GA, 2007GA}, evolutionary programming (EP) \cite{1999EP, 2004EP, 2011EP, 2016EP}, particle swarm optimization (PSO) \cite{2003PSO, 2007PSO, 2010PSO, 2015PSO}, $\lambda$-iterative technique \cite{2014lambda}, nonlinear programming (NLP) \cite{2010NLP}, semi-definite programming (SDP) \cite{2012SDP}, mixed integer quadratic programming (MIQP) \cite{2007MIQP, 2014BigM, 2016Distance}, branch and bound (B\&B) \cite{2015BB}, etc. Although various optimization methods have been used to overcome the discontinuous solution space of SED-POZ, it still can not be easily implemented for the solution of DED-POZ. When POZ is taken into account, complexity of the DED problem will increase significantly since DED dispatches over a scheduled time horizon instead of one period. Moreover, additional constraints such as ramp rate constraints and spinning reserve constraints make the problem more complicated and can not be tackled easily.

More recent works for DED-POZ have been around heuristic methods, such as particle swarm optimization (PSO) \cite{2004meetingPOZ, 2009yuan, 2012PSO}, imperialist competitive algorithm (ICA) \cite{2012Imperialist}, harmony search (HS) \cite{2014HS}, improved bees algorithm (IBA) \cite{2015Improved}, etc. However, heuristic techniques are quite sensitive to various parameter settings and solution may be different at each trial due to the intrinsic stochastic characteristic of heuristics. They do not provide an optimality gap so you have no clue how well of a solution you have obtained. Hybrid methods which combine several heuristic techniques or deterministic approach such as bacterial foraging PSO-DE algorithm (BPSO-DE) \cite{2012PSODE} and hybrid EP-PSO-SQP algorithm \cite{2008hybrid} tend to be more efficient than the individual method. However, they still have the intrinsic drawback of the stochastic search method we mentioned above. Unlike heuristics, deterministic algorithms can solve to a robust result due to the solid mathematical foundations and the availability of the powerful software tools. In \cite{2014optimality}, an efficient real-time approach based on optimality condition decomposition (OCD) technique is proposed to solve the DED-POZ. By using the reformulation and OCD technique, the problem could be decomposed to several simpler sub-problems and then the CPU-time can be reduced significantly, but the spinning reserve constraints are not considered. In \cite{2015MIQP}, the DED-POZ formulates an MIQP model which can be solved by a mixed integer programming (MIP) solver immediately. Nevertheless, the complicated transmission losses are not included.

As the significant progresses of mathematical programming theory and the improvements of the MIP solvers, the MIP technique has become a popular alternative in the optimal operation of electric power system. In \cite{2007MIQP, 2014BigM, 2016Distance, 2015MIQP}, various MIQP formulations which are capable to solve to global optimality directly are presented for SED-POZ and DED-POZ, while the transmission losses are neglected. Although the solution of MIQP with commercial software such as CPLEX or GUROBI have significantly improved in recent years, it still does not a very good choice for DED-POZ on account of its non-linear characteristic. Comparing with the MIQP, mixed integer linear programming (MILP) is more developed because of the vastly superior warm-start capabilities of the simplex method \cite{2012Jabr}. Besides, from the perspective of engineering, an exact optimal solution of the problem is not always necessary and faster approximations have more value. Consequently, a natural alternative for solving DED-POZ is to reformulate the problem, yielding a reliable MILP formulation which can be solved efficiently by an MIP solver.

To find such an MILP formulation for DED-POZ, auxiliary variables are introduced not only for the constraints but also the objective function in our work. By using the perspective cut reformulation technique, a tight MILP formulation, denoted by MILP-1, which can be solved to global optimality within a preset tolerance via an MIP solver is obtained. When transmission losses are considered, the complicated quadratic terms in the transmission losses are replaced by their first order Taylor expansions. Then an MILP formulation for DED considering transmission losses, denoted by MILP-2, is formed. Based on MILP-1 and MILP-2, we propose an MILP-iteration algorithm (MILP-IA) to solve the DED problem. Since the previous studies mainly focus on the optimality of the solution and did not discuss the feasibility in detail for the DED-POZ. In our work, not only the optimality but also the feasibility are discussed. Simulation results show that the proposed MILP-1 and MILP-IA can solve to competitive solutions in a short time.

The rest of this paper is organized as follows. The mathematical formulation for the DED-POZ is described in Section \ref{Mathematical}. The reformulation and MILP-IA are presented in Section \ref{Reformulation}. Simulation results on several test systems are given in Section \ref{Simulation}. Finally, conclusions are drawn in Section \ref{Conclusion}.

\section{Mathematical formulation for the DED-POZ}
\label{Mathematical}
The objective of DED-POZ is to minimize the total generation cost over the scheduled time horizon
\begin{equation}\label{objective}
\min~\sum\limits_{t=1}^T\sum\limits_{i=1}^N(\alpha_i+\beta_iP_{i,t}+\gamma_iP_{i,t}^2)
\end{equation}
where $P_{i,t}$ is the power output of unit $i$ in period $t$; $N$ is the total number of units; $T$ is the total number of periods; $\alpha_i$, $\beta_i$ and $\gamma_i$ are the cost coefficients of unit $i$.

The minimized DED-POZ should be subjected to the constraints as follows.

1) Power balance equations
\begin{equation}\label{balance equation}
\sum\limits_{i=1}^N P_{i,t}=D_{t}+P_{t}^{loss}, ~~~\forall~t
\end{equation}
where $D_{t}$ is the load demand in period $t$; $P_{t}^{loss}$ is the transmission loss in period $t$, which can be calculated based on Kron's loss formula as follow \cite{1999Saadat}:
\begin{equation}\label{loss}
P_{t}^{loss}=B_{00}+B_{0}^{T}P_{t}+P_{t}^{T}BP_{t}, ~~~\forall~t
\end{equation}
where $B_{00}$, $B_0$, $B$ are B-coefficients; $P_{t}=[P_{1,t},P_{2,t},...,P_{N,t}]^T$ is the power output vector in period $t$.

2) Power generation limits
\begin{equation}\label{generation limits}
P_i^{\min}\le{P_{i,t}}\le{P_i^{\max}}, ~~~\forall~i,t
\end{equation}
where $P_i^{\min}$ and $P_i^{\max}$ are the minimum and maximum power outputs of unit $i$.

3) Ramp rate limits
\begin{equation}\label{ramp limits}
RD_{i}\le{P_{i,t}-P_{i,t-1}}\le{RU_{i}}, ~~~\forall~i,t
\end{equation}
where $RD_{i}$ and $RU_{i}$ are the ramp-down and ramp-up rates of unit $i$.

4) Prohibited operating zones limits

The prohibited operating zones of each unit can be characterised by the disjoint operating regions as shown below
\begin{equation}\label{POZ Limits}
\left\{
\begin{array}{l}
P_{i1}^{\min}\le{P_{i,t}}\le{P_{i1}^{\max}} ~~\textrm{or}  \\
P_{ij}^{\min}\le{P_{i,t}}\le{P_{ij}^{\max}} ~~\textrm{or}  \\
P_{i{n_i+1}}^{\min}\le{P_{i,t}}\le{P_{i{n_i+1}}^{\max}},~j=2, ..., n_i,~\forall~i,t
\end{array}
\right.
\end{equation}
where $P_{i1}^{\min}=P_i^{\min}$, $P_{i{n_i+1}}^{\max}=P_i^{\max}$, $n_i$ is the number of prohibited operating zones of unit $i$.

5) Spinning reserve constraints
\begin{equation}\label{reserve requirements}
\left\{
\begin{array}{l}
SR_{i,t}\le \min {\{P_i^{\max}-P_{i,t}, RU_{i}\}}, ~~~\forall~i,t  \\
\sum\limits_{i=1}^N SR_{i,t}\geq R_{t}, ~~~\forall~t
\end{array}
\right.
\end{equation}
where $SR_{i,t}$ is the spinning reserve provided by unit $i$ in period $t$ and $R_{t}$ is the system spinning reserve requirement in period $t$.

\section{Reformulation and solution for the DED-POZ}
\label{Reformulation}
Owing to the disjoint operating zones, the classical mathematical programming methods are not suitable for DED-POZ any more. By introducing some auxiliary variables $P_{i,t}^{j}$ and binary variables $u_{i,t}^{j}(j=1,...,n_i+1)$, the POZ constraint (\ref{POZ Limits}) can be equivalent to the following expressions
\begin{equation}\label{new POZ Limits 1}
\left\{
\begin{array}{l}
u_{i,t}^{1}P_{i1}^{\min}\le{P_{i,t}^{1}}\le{u_{i,t}^{1}P_{i1}^{\max}}   \\
u_{i,t}^{j}P_{ij}^{\min}\le{P_{i,t}^{j}}\le{u_{i,t}^{j}P_{ij}^{\max}}    \\
u_{i,t}^{n_i+1}P_{i{n_i+1}}^{\min}\le{P_{i,t}^{n_i+1}}\le{u_{i,t}^{n_i+1}P_{i{n_i+1}}^{\max}},~j=2, ..., n_i,~\forall~i,t  \\
\end{array}
\right.
\end{equation}
\begin{equation}\label{new POZ Limits 2}
\sum\limits_{j=1}^{n_i+1}P_{i,t}^{j}=P_{i,t},~\forall~i,t
\end{equation}

\begin{equation}\label{new POZ Limits 3}
\left\{
\begin{array}{l}
\sum\limits_{j=1}^{n_i+1}u_{i,t}^{j}=1  \\
u_{i,t}^{j}\in\{0,1\},~\forall~i,t.   \\
\end{array}
\right.
\end{equation}

Consequently, the DED-POZ can be formulated as an MIQP formulation when the transmission losses are not included:
\begin{equation}\label{MIQP}
\begin{array}{l}
\min~\sum\limits_{t=1}^T\sum\limits_{i=1}^N(\alpha_i+\beta_iP_{i,t}+\gamma_iP_{i,t}^2) \\
s.t.~(\ref{balance equation}),(\ref{generation limits}),(\ref{ramp limits}),(\ref{reserve requirements}),(\ref{new POZ Limits 1}),
(\ref{new POZ Limits 2}),(\ref{new POZ Limits 3}).
\end{array}
\end{equation}
The DED-POZ, as formulated in (\ref{MIQP}), can be solved via MIQP directly. However, solution via MILP tends to be more efficient since the warm start capabilities of the simplex method available in MILP solver are vastly superior in comparison with the interior-point method in MIQP solver. Moreover, from the perspective of engineering, an exact optimal solution of the problem is not always necessary and faster approximations have more value. As a result, we present two reliable MILP formulations for solving the DED-POZ in the next subsections.

\subsection{ Reformulation for the DED-POZ}

According to $u_{i,t}^{j}P_{ij}^{\min}\le{P_{i,t}^{j}}\le{u_{i,t}^{j}P_{ij}^{\max}}$ and $u_{i,t}^{j}\in\{0,1\}$, we have $P_{i,t}^{j}\in\{0\}\cup[P_{ij}^{\min},P_{ij}^{\max}]$, which indicates that $P_{i,t}^{j}$ is a semi-continuous variable. With the help of (\ref{new POZ Limits 2}), the objective function in (\ref{MIQP}) can be converted into the sum of quadratic functions on the semi-continuous variables space
\begin{eqnarray}\label{objective_1}
\min~\sum\limits_{t=1}^T\sum\limits_{i=1}^N(\alpha_i+\sum\limits_{j=1}^{n_i+1}(\beta_iP_{i,t}^{j}+\gamma_i(P_{i,t}^{j})^2)).
\end{eqnarray}
Perspective cut reformulation technique \cite{2006Frangioni, 2009Frangioni} proposed by Frangioni et al. can thus be used for constructing a tight MILP formulation of the problem. Introducing some auxiliary variables $z_{i,t}^{j}(j=1,...,n_i+1)$, then (\ref{objective_1}) can be approximated by the following perspective cuts
\begin{eqnarray}\label{objective_2}
\min~\sum\limits_{t=1}^T\sum\limits_{i=1}^N(\alpha_i+\sum\limits_{j=1}^{n_i+1} z_{i,t}^{j}) ~~~~~~~~~~~~~~~~~~~~~~~~~~~~~~~~~\\
\label{objective_pcz}
s.t.~z_{i,t}^{j}\geq (2\gamma_i P_{i}^{j(l)}+\beta_i)P_{i,t}^{j}-\gamma_i(P_{i}^{j(l)})^2u_{i,t}^{j},~\forall~i,t,j
\end{eqnarray}
where $P_{i}^{j(l)}=P_{ij}^{\min}+l(P_{ij}^{\max}-P_{ij}^{\min})/L$$(l = 0,1,...,L)$ and $L$ is a given parameter.

Then a tight MILP approximation, denoted as MILP-1, for DED-POZ without considering transmission losses can be formed:
\begin{equation} \label{PC-MILP}
\begin{array}{l}
\min~\sum\limits_{t=1}^T\sum\limits_{i=1}^N(\alpha_i+\sum\limits_{j=1}^{n_i+1} z_{i,t}^{j}) \\
s.t.~(\ref{balance equation}),(\ref{generation limits}),(\ref{ramp limits}),(\ref{reserve requirements}),(\ref{new POZ Limits 1}),
(\ref{new POZ Limits 2}),(\ref{new POZ Limits 3}),(\ref{objective_pcz}).
\end{array}
\end{equation}

\subsection{Linearization for the transmission losses}

When the transmission losses are taken into account for the DED-POZ, it makes the optimization more difficult because of its complicated nonlinearity. Note that, the nonlinearity arises from the third term of (\ref{loss}), i.e., $P_{t}^{T}BP_{t}$. Replacing $P_{t}^{T}BP_{t}$ with auxiliary variable $c_t$, the transmission loss (\ref{loss}) can be rewritten as
\begin{eqnarray}\label{loss1}
P_{t}^{loss}=B_{00}+B_{0}^{T}P_{t}+c_{t}, ~~~\forall~t     \\
\label{c t0}
c_{t} \geq P_{t}^{T}BP_{t}, ~~~~~~~~~\forall t.~~~~~~
\end{eqnarray}

As we can see, (\ref{loss1}) is a linear constraint which can be addressed easily while (\ref{c t0}) is a complicated quadratic constraint which is hard to tackle. A natural way to conquer this difficulty is to find an efficient linear approximation instead of the quadratic one. Then the first order Taylor expansion is employed for the $P_{t}^{T}BP_{t}$. As a result, the (\ref{c t0}) can be replaced by
\begin{equation}\label{c_t}
c_{t} \geq 2P_{t}^{(k)T}BP_{t}-P_{t}^{(k)T}BP_{t}^{(k)}, ~~~\forall~t
\end{equation}
where $P_{t}^{(k)}$ is taken to be a constant vector corresponding to $P_{t}$.

When the transmission losses are considered, based on MILP-1 and the above approximations, the DED-POZ can be formulated as the following MILP formulation, denoted as MILP-2
\begin{equation}\label{PCL}
\begin{array}{l}
\min~\sum\limits_{t=1}^T\sum\limits_{i=1}^N(\alpha_i+\sum\limits_{j=1}^{n_i+1} z_{i,t}^{j}) \\
s.t.~(\ref{balance equation}),(\ref{generation limits}),(\ref{ramp limits}),(\ref{reserve requirements}),(\ref{new POZ Limits 1}),
(\ref{new POZ Limits 2}),          \\
~~~~~(\ref{new POZ Limits 3}),(\ref{objective_pcz}),
(\ref{loss1}),(\ref{c_t}).
\end{array}
\end{equation}

It is well know that, when the vector $P_{t}^{(k)}$ we employ in (\ref{c_t}) is sufficiently close to the optimal solution, a high efficiency MILP approximation for the primal problem can be obtained. Then MILP-2 can provide a near-optimal solution. To exploit such an efficient MILP approximation and solve to a reliable solution, based on MILP-1 and MILP-2 formulation, a straightforward MILP-iteration algorithm, denoted as MILP-IA, is proposed for DED-POZ when transmission losses are considered. In the following, the details of MILP-IA are given.

\textbf{Initialization Step:} Choose a scalar $\epsilon > 0$ and a maximum number of iteration $iter_{max}$ to be used for terminating the algorithm and let $iter=1$, $k=3$.

\textbf{Step 1}: Solve MILP-1 (\ref{PC-MILP}) to obtain an optimal solution denoted by $P^{(1)}$, where
$P^{(1)}=[P_{1}^{(1)};P_{2}^{(1)};...;P_{T}^{(1)}]$.

\textbf{Step 2}: Solve MILP-2 (\ref{PCL}) where the linear approximation (\ref{c_t}) is taken at $P_{t}^{(1)}$, to obtain an optimal solution denoted by $P^{(2)}$.

\textbf{Step 3}: Solve MILP-2 (\ref{PCL}) where the linear approximations (\ref{c_t}) is taken at $(P_{t}^{(k-2)}+P_{t}^{(k-1)})/2$ , to obtain an optimal solution denoted by $P^{(k)}$.

\textbf{Step 4}:  Calculate the violation of power balance $E_t$, where
\begin{equation}\label{error_t}
E_t=|\sum\limits_{i=1}^N P_{i,t}^{(k)}-D_{t}-P_{t}^{loss(k)}|,~\forall~t
\end{equation}
and the $P_{t}^{loss(k)}$ is calculated according to (\ref{loss}).

\textbf{Step 5}:  When
\begin{equation}\label{stop}
E_t < \epsilon, ~\forall~t
\end{equation}
or when $iter=iter_{max}$, the procedure is terminated.

\textbf{Step 6}:  $iter=iter+1$, $k = k+1$. Go to step 3.

Because the transmission loss at each period in a DED problem is small compared to the corresponding load demand, when transmission losses are ignored, after optimizing MILP-1 (\ref{PC-MILP}), an initial approximate optimal solution can be obtained in step 1. Since transmission losses are ignored, the solution obtained in step 1 is ``relative small" in some ways. When transmission losses are included, by solving MILP-2 (\ref{PCL}) where the linear approximation (\ref{c_t}) is taken at such a ``relative small" solution, a ``relative large" solution can be got in step 2. Then average of the former two solutions is used in step 3 to balance the power balance equation. Repeat from step 3 to 6 until all the violations of power balances are smaller than the preset value or the maximum number of iteration is reached. Generally, several iterations are enough in our algorithm since the transmission loss is small and $P_{t}^{T}BP_{t}$ is only a portion of the transmission loss.

\section{Simulation results}
\label{Simulation}
To assess the efficiency of the proposed MILP-1 formulation and MILP-IA for DED-POZ, two cases, ignoring the transmission losses and considering the transmission losses, are simulated in our numerical experiments. The formulation and algorithm are coded with Matlab and optimized by using CPLEX 12.6.2 \cite{2015CPLEX}. The machine for all runs is an Intel Core 2.5 GHz Dell-notebook with 8 GB of RAM.

\subsection{Ignoring the transmission losses}

In this subsection, a set of different sizes test systems with units ranging from 6 to 180 over a scheduled time horizon of 24 h are adopted for testing the effectiveness of MILP-1, where the transmission losses are not considered. The 6-unit system data is taken from \cite{2004meetingPOZ}. The 1-h spinning reserve requirement is 5\% of the load demand. Other test systems with 30, 60, 120, and 180 units can be obtained by duplicating the 6-unit system 5, 10, 20, and 30 times, respectively. In our MILP-1 formulation, $L$ is set to 4.

\renewcommand{\multirowsetup}{\centering}
\begin{table}[!ht]
\caption{\label{Table24} Comparison of the two formulations}
\centering
\resizebox{0.55\textwidth}{!}{
\begin{tabular}{c|cc|cc}
\hline
\multirow{2}{*}{N} & \multicolumn{2}{c}{MIQP} &  \multicolumn{2}{|c}{MILP-1}  \\
\cline{2-5}
& Cost(\$)& Time(s)  & Cost(\$)& Time(s)    \\
\hline
6&         310500&           1.23&       310506&     0.22  \\
30&        1552531&         45.51&      1552541&     0.71  \\
60&        3105094&        154.95&      3105088&     1.31  \\
120&       6210747&        300.00&      6210174&     2.93  \\
180&       9315969&        300.00&      9315270&     4.84  \\
\hline
\end{tabular}}
\end{table}

We directly solve the MIQP (\ref{MIQP}) and MILP-1 (\ref{PC-MILP}) formulations using CPLEX to 0.01\% optimality in a time limit of 300 s, and the results are compared in Table \ref{Table24}. As we can see in Table \ref{Table24}, CPLEX can solve MILP-1 faster than it solves MIQP, with much fewer computation time. For MILP-1, solutions for all systems can be gained in several seconds, while for MIQP, only the 6-, 30- and 60- unit systems can be solved to the preset tolerance within 300 s. On the other hand, although the total generation costs for the 10- and 30- unit systems obtained by solving MILP-1 are larger than the ones obtained by solving MIQP, but the reduced costs are very small. Furthermore, for the 60-, 120- and 180- unit systems, MILP-1 can exploit more lower total generation costs within several seconds.

The outputs obtained by solving MILP-1 for the 6-unit system are given in Table \ref{Tableoutputs10unit} for verification.

\begin{table}[!ht]
\caption{\label{Tableoutputs10unit} Outputs (MW) for the 6-unit system }
\centering
\resizebox{0.6\textwidth}{!}{
\begin{tabular}{c|cccccc}
\hline
\multicolumn{1}{c|}{\raisebox{0.1ex}[0cm][0cm]{$t$}}& unit $1$&  unit $2$&  unit $3$&  unit $4$&  unit $5$& unit $6$  \\
\hline
1& 383.75& 121.25& 210.00&  76.25& 113.75 & 50.00 \\
\hline
 2& 380.00& 121.25& 208.25&  68.75& 113.75 & 50.00 \\
\hline
 3& 380.00& 121.25& 205.00&  68.75& 110.00 & 50.00 \\
\hline
 4& 380.00& 116.25& 205.00&  68.75& 110.00 & 50.00 \\
\hline
 5& 380.00& 121.25& 205.00&  68.75& 110.00 & 50.00 \\
\hline
 6& 391.75& 121.25& 210.00&  76.25& 113.75 & 50.00 \\
\hline
 7& 395.00& 128.75& 210.00&  80.00& 125.25 & 50.00 \\
\hline
 8& 395.00& 139.25& 210.00&  92.50& 136.25 & 50.00 \\
\hline
 9& 425.00& 140.00& 247.50& 104.12& 150.00 & 59.38 \\
\hline
10& 425.00& 160.00& 247.50& 107.50& 150.62 & 59.38 \\
\hline
11& 425.00& 165.00& 262.50& 120.00& 156.63 & 71.88 \\
\hline
12& 440.00& 165.00& 262.50& 123.75& 168.75 & 75.00 \\
\hline
13& 425.00& 165.00& 251.88& 120.00& 156.25 & 71.88 \\
\hline
14& 455.00& 166.00& 262.50& 123.75& 168.75 & 75.00 \\
\hline
15& 455.00& 168.00& 262.50& 123.75& 168.75 & 85.00 \\
\hline
16& 455.00& 165.00& 262.50& 123.75& 168.75 & 75.00 \\
\hline
17& 429.75& 165.00& 262.50& 120.00& 168.75 & 75.00 \\
\hline
18& 425.00& 165.00& 262.50& 120.00& 157.63 & 71.88 \\
\hline
19& 425.00& 160.00& 247.50& 107.50& 156.25 & 62.75 \\
\hline
20& 425.00& 140.00& 240.00&  97.50& 139.50 & 50.00 \\
\hline
21& 395.00& 139.25& 210.00&  92.50& 136.25 & 50.00 \\
\hline
22& 395.00& 128.75& 210.00&  79.00& 121.25 & 50.00 \\
\hline
23& 395.00& 128.75& 210.00&  76.25& 115.00 & 50.00 \\
\hline
24& 388.75& 121.25& 210.00&  76.25& 113.75 & 50.00 \\
\hline
\end{tabular}}
\end{table}

\subsection{Considering the transmission losses}

In this subsection, a 6-unit system and a 15-unit system, which have been widely used for testing \cite{2004meetingPOZ, 2009yuan, 2015Improved} are taken into account.
Although the systems have only 6 units and 15-units, they are dynamic problem with 288 and 750 variables coded in a solution, respectively. The characteristics of the thermal units and load demands are taken from \cite{2004meetingPOZ}. Owing to the limits of space, the loss coefficients with the 100-MVA base capacity are not listed here. One can refer to \cite{2003PSO}. For fair comparison, the spinning reserve requirement is 5\% of the load demand. In our MILP-IA, parameters $\epsilon$, $iter_{max}$ and $L$ are set to 0.1, 5 and 4, respectively. Meanwhile, in MILP-IA, we directly solve the MILP-1 and MILP-2 formulations using CPLEX to 0.01\% optimality.

Since different computers are used, the run time in different methods may not be directly comparable. In order to have a fair comparison regarding the computational effort, the CPU chip frequency from the used computer is used to convert the CPU times obtained from different methods into a common base \cite{2008scaletime}. The scaled CPU time of an algorithm can be computed as follow:
\begin{equation}\label{P2pwan}
Scaled~CPU~time = \frac{Given~CPU~speed}{Base~CPU~speed} Given~CPU~time.
\end{equation}
For fair comparison, the scaled CPU time is used in this subsection and the base CPU speed is 2.0 GHz.

\subsubsection{6-unit system considering the transmission losses}

In this subsection, a 6-unit system considering the transmission losses is adopted for simulation. In this system, all the units have POZ constraints. The results obtained by MILP-IA and IBA \cite{2015Improved} are shown in Table \ref{Table6unitcompare}. We can see that our MILP-IA can save much more execution time than IBA. Although IBA can obtain much lower total generation cost, but we should note that, not only the optimality but also feasibility should be taken into account for evaluating the quality of the solution.

\begin{table}[!ht]
\caption{\label{Table6unitcompare} Results of MILP-AI and IBA for the 6-unit system}
\centering
\resizebox{0.36\textwidth}{!}{
\begin{tabular}{c|c|c}
\hline
Method&       Cost(\$)&   Time(s)   \\
\hline
IBA \cite{2015Improved}   &    313472  &  13.2    \\
\hline
 MILP-IA &       315169      &   0.84  \\
\hline
\end{tabular}}
\end{table}

In DED, the constraints such as power generation limits, ramp rate limits, POZ limits and spinning reserve constraints are linear constraints which can be easily satisfied in most algorithms. But the power balance equations are quadratic equality constraints which are hard to meet entirely. In fact, thanks to the MIP solver, the feasibility of those linear constraints can be easily guaranteed in MILP-IA. Therefore, our discussion with respect to the feasibility validation mainly focus on the violations of power balances.

Since the losses coefficients are specified in per unit on a 100-MVA base in our test systems, then the real power loss is given by \cite{1999Saadat}:
\begin{equation}\label{losspu}
P_{t(pu)}^{loss}=B_{00}+B_{0}^{T}P_{t(pu)}+P_{t(pu)}^{T}BP_{t(pu)}, ~~~\forall~t.
\end{equation}
In the cost function $P_{i,t}$ is expressed in MW. Therefore, the real power loss in terms of MW generation is \cite{1999Saadat}:
 \begin{equation}\label{lossMW}
P_{t}^{loss}=[B_{00}+B_{0}^{T}(\frac{P_{t}}{100})+(\frac{P_{t}}{100})^{T}B(\frac{P_{t}}{100})]\times100 ~\textrm{MW}, ~~~\forall~t.
\end{equation}
Hence, by using (\ref{lossMW}) and (\ref{error_t}), the violations of power balances can be calculated easily.

The comparison results with respect to the violations for the 6-unit system between IBA and MILP-IA are shown in Table \ref{Table6unitviolation}.

\begin{table}[!ht]
\caption{\label{Table6unitviolation} Violations (MW) for the 6-unit system}
\centering
\resizebox{0.75\textwidth}{!}{
\begin{tabular}{c|cccccccc}
\hline
$t$ &1&2&3&4&5&6&7&8    \\
\hline
IBA \cite{2015Improved}& 5.5392&  5.5574&  5.5495&  5.5367&  5.5410&  5.5393&	5.5466&  5.5412 \\
MILP-IA& 0.0004&  0.0008&  0.0008&  0.0005&  0.0008&  0.0011&  0.0018&  0.0161\\
\hline
$t$&9&10&11&12&13&14&15&16    \\
\hline
IBA \cite{2015Improved}&  5.5551&  5.5044&  5.5390&  5.5456& 5.5389&  5.5412&  5.5312&  5.5424 \\
MILP-IA&  0.0065&  0.0035&  0.0037&  0.0101&  0.0010&  0.0095&  0.0021&  0.0095  \\
\hline
$t$&17&18&19&20&21&22&23&24    \\
\hline
IBA \cite{2015Improved}&  5.5450&  5.5519& 5.5418&  5.5216&  5.5598&  5.5392&  5.5323&  5.5255  \\
MILP-IA&  0.0101&  0.0039&  0.0040&  0.0040&  0.0161&  0.0011&  0.0012&  0.0004  \\
\hline
\end{tabular}}
\end{table}

As we can see in Table \ref{Table6unitviolation}, almost all the violations of power balances for IBA are more than 5 MW, which makes the solution infeasible. Whereas, all the violations of power balances for our MILP-IA are less than 0.1 MW. Actually, when the DED-POZ is solved by using MILP-IA, the algorithm terminates at the second iteration. Moreover, in Table \ref{Table6unitviolation}, most of the violations for MILP-IA are much smaller than 0.1 MW and the total violation is only 0.1090 MW, which shows that MILP-IA can solve to an acceptable near-optimal solution for the 6-unit system within 1 s.

The outputs and transmission losses for the 6-unit system obtained by MILP-IA for DED-POZ are given in Table \ref{Table6unitOutputs} for verification.

\begin{table}[!ht]
\caption{\label{Table6unitOutputs} Outputs and losses (MW) for the 6-unit system}
\centering
\resizebox{0.6\textwidth}{!}{
\begin{tabular}{c|ccccccc}
\hline
\multicolumn{1}{c|}{\raisebox{0.1ex}[0cm][0cm]{$t$}}& unit $1$&  unit $2$&  unit $3$&  unit $4$&  unit $5$&  unit $6$&  loss \\
\hline
1& 393.18& 121.25& 210.00&  80.00& 113.75&  50.00&  13.19  \\
\hline
 2& 380.00& 121.25& 209.96&  80.00& 113.75&  50.00&  12.96  \\
\hline
 3& 380.00& 121.25& 205.00&  80.00& 111.58&  50.00&  12.83  \\
\hline
 4& 380.00& 121.25& 205.00&  76.54& 110.00&  50.00&  12.80  \\
\hline
 5& 380.00& 121.25& 205.00&  80.00& 111.58&  50.00&  12.83  \\
\hline
 6& 391.16& 121.25& 210.00&  90.00& 113.75&  50.00&  13.16  \\
\hline
 7& 395.00& 133.82& 210.00&  92.50& 121.25&  50.00&  13.58  \\
\hline
 8& 395.00& 136.25& 240.00&  94.69& 121.25&  50.00&  14.21  \\
\hline
 9& 425.00& 160.00& 241.15& 110.00& 140.00&  65.62&  15.78  \\
\hline
10& 425.00& 160.00& 245.51& 120.00& 150.00&  65.62&  16.14  \\
\hline
11& 427.96& 165.00& 262.50& 131.25& 156.25&  75.00&  16.97  \\
\hline
12& 452.60& 165.00& 262.50& 131.25& 156.25&  85.00&  17.61  \\
\hline
13& 425.00& 165.00& 254.23& 131.25& 156.25&  75.00&  16.73  \\
\hline
14& 455.00& 171.29& 262.50& 138.75& 156.25&  85.00&  17.80  \\
\hline
15& 455.00& 175.00& 262.50& 138.75& 164.85&  85.00&  18.11  \\
\hline
16& 455.00& 170.27& 262.50& 138.75& 156.25&  85.00&  17.78  \\
\hline
17& 438.31& 165.00& 262.50& 131.25& 156.25&  85.00&  17.32  \\
\hline
18& 428.98& 165.00& 262.50& 131.25& 156.25&  75.00&  16.99  \\
\hline
19& 425.00& 160.00& 247.50& 123.75& 150.00&  68.98&  16.23  \\
\hline
20& 424.11& 140.00& 240.00& 107.50& 136.25&  59.37&  15.23  \\
\hline
21& 395.00& 136.25& 240.00&  94.69& 121.25&  50.00&  14.21  \\
\hline
22& 395.00& 128.75& 210.00&  92.50& 121.25&  50.00&  13.50  \\
\hline
23& 395.00& 128.75& 210.00&  90.00& 114.61&  50.00&  13.36  \\
\hline
24& 395.00& 124.52& 210.00&  80.00& 113.75&  50.00&  13.27  \\
\hline
\end{tabular}}
\end{table}

\subsubsection{15-unit system considering the transmission losses}

In this subsection, a 15-unit system considering the transmission losses is adopted for simulation. In this system only 4 units have POZ constraints. The results obtained by MILP-IA and other two PSO methods \cite{2009yuan} are shown in Table \ref{Table15unitcompare}.

\begin{table}[!ht]
\caption{\label{Table15unitcompare} Results of MILP-AI and other methods for the 15-unit system }
\centering
\resizebox{0.32\textwidth}{!}{
\begin{tabular}{c|c|c}
\hline
Method&       Cost(\$)&     Time(s)     \\
\hline
PSO \cite{2009yuan}  &  761774 &     6.00\\
EPSO \cite{2009yuan} &  759410 &     7.50\\
MILP-IA  &             759176 &     1.03     \\
\hline
\end{tabular}}
\end{table}

\begin{figure}[!ht]
\centering\includegraphics[width=2.8in]{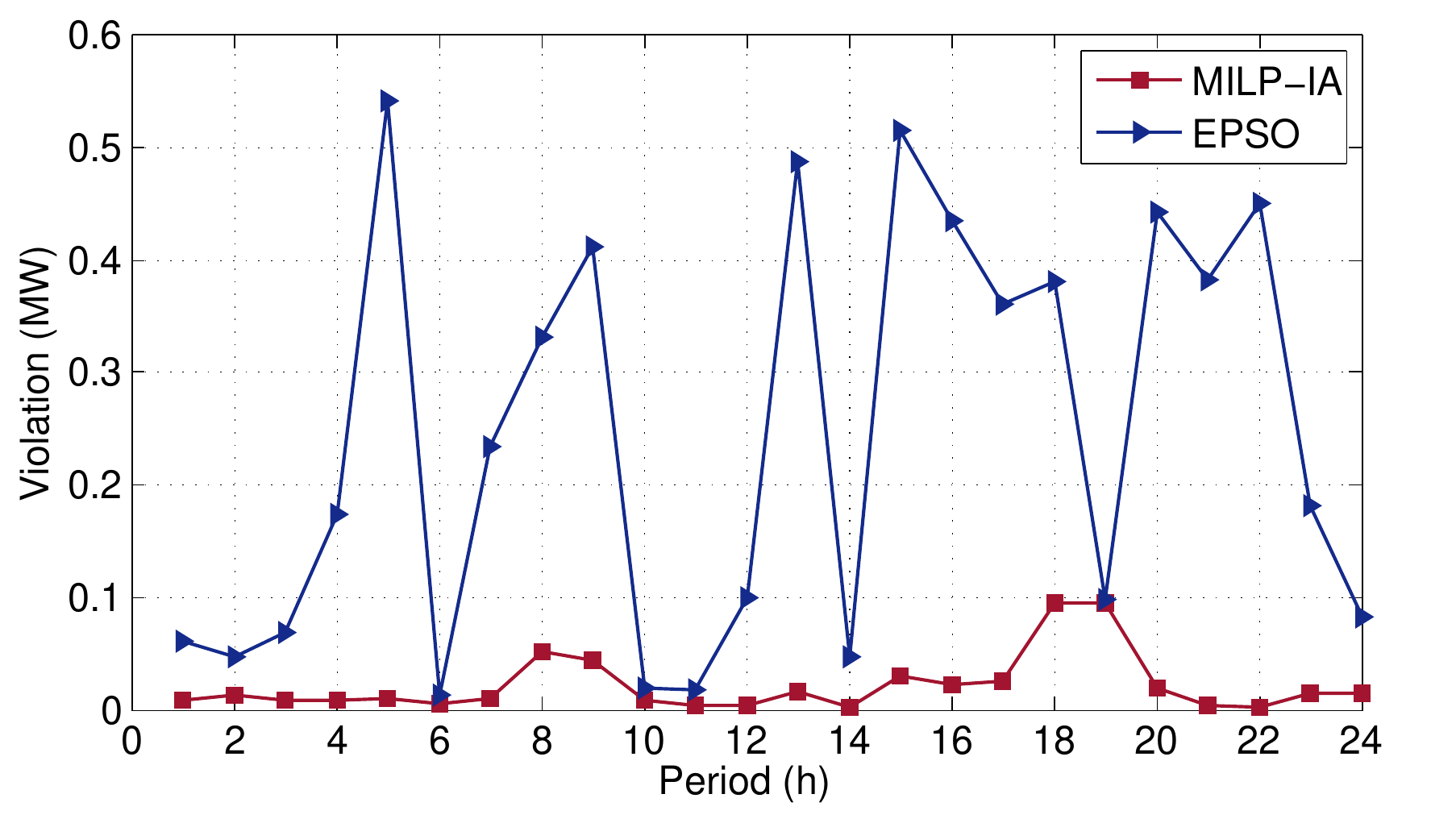}
\caption{Violations of power balances for the 15-unit system}\label{epso}
\end{figure}

As we can see, our MILP-IA can solve to a lower total generation cost in a faster speed than other two PSO methods. As stated before, lower cost does not mean a good solution. The violations of power balances must be checked. In Table \ref{Table15unitdifferent}, only the violations for EPSO and MILP-IA are calculated, since the solutions for PSO are not available in \cite{2009yuan}. Figure \ref{epso} depicts the details as well.

\begin{table}[!ht]
\caption{\label{Table15unitdifferent} Violations (MW) for the 15-unit system }
\centering
\resizebox{0.75\textwidth}{!}{
\begin{tabular}{c|cccccccc}
\hline
$t$ &1&2&3&4&5&6&7&8   \\
\hline
EPSO\cite{2009yuan}&  0.0612&  0.0476&  0.0684&  0.1747&  0.5416&  0.0127&  0.2344&  0.3323 \\
MILP-IA&   0.0090&  0.0135&  0.0083&  0.0090&  0.0098&  0.0047&  0.0095&  0.0513 \\
\hline
$t$ &9&10&11&12&13&14&15&16   \\
\hline
EPSO\cite{2009yuan}&  0.4119&  0.0196&  0.0185&  0.0998&  0.4874&  0.0473&  0.5161&  0.4348\\
MILP-IA&    0.0448&  0.0086&  0.0039&  0.0031&  0.0163&  0.0021&  0.0301&  0.0224 \\
\hline
$t$&17&18&19&20&21&22&23&24    \\
\hline
EPSO\cite{2009yuan}&  0.3617&  0.3817&  0.0976&  0.4436&  0.3827&  0.4506&  0.1822&  0.0826 \\
MILP-IA&     0.0258&  0.0949&  0.0946&  0.0198&  0.0038&  0.0019&  0.0151&  0.0152\\
\hline
\end{tabular}}
\end{table}

We can see from Table \ref{Table15unitdifferent} and Figure \ref{epso}, all the violations for MILP-IA are smaller than EPSO and most of them are much smaller than EPSO. In EPSO, the maximum violation is 0.5416 MW. But in our algorithm, the maximum violation is 0.0949 MW and the total violation is only 0.5175 MW. Actually, when DED-POZ is solved by our MILP-IA, the algorithm terminates at the second iteration. In other words, the process terminates because all the violations of power balances are less than 0.1 MW. From the simulation results we can conclude that our MILP-IA can solve to a better solution in a short time for DED-POZ.

The outputs and transmission losses for the 15-unit system obtained by MILP-IA are given in Table \ref{Table15unitoutputs} for verification.

\begin{table}[!ht]
\caption{\label{Table15unitoutputs} Outputs and losses (MW) for 15-unit system}
\resizebox{1\textwidth}{!}{
\begin{tabular}{c|cccccccccccccccc}
\hline
\multicolumn{1}{c|}{\raisebox{0.1ex}[0cm][0cm]{$t$}}& unit $1$&  unit $2$&  unit $3$&  unit $4$&  unit $5$&  unit $6$&  unit $7$&  unit $8$& unit $9$&  unit $10$&  unit $11$&  unit $12$&  unit $13$&  unit $14$&  unit $15$&  loss \\
\hline
 1& 369.25& 295.00& 130.00& 130.00& 150.00& 455.63& 465.00&  60.00&  25.00&  25.00&  42.50&  53.13&  25.00&  15.00&  15.00&  19.51 \\
\hline
 2& 352.34& 295.00& 130.00& 130.00& 150.00& 455.00& 465.00&  60.00&  25.00&  25.00&  42.50&  49.37&  25.00&  15.00&  15.00&  19.23 \\
\hline
 3& 359.72& 295.00& 130.00& 130.00& 150.00& 455.00& 465.00&  60.00&  25.00&  25.00&  42.50&  53.13&  25.00&  15.00&  15.00&  19.35 \\
\hline
 4& 369.25& 295.00& 130.00& 130.00& 150.00& 455.63& 465.00&  60.00&  25.00&  25.00&  42.50&  53.13&  25.00&  15.00&  15.00&  19.51 \\
\hline
 5& 416.88& 305.00& 130.00& 130.00& 150.00& 460.00& 465.00&  60.00&  25.00&  25.00&  42.50&  54.15&  25.00&  15.00&  15.00&  20.53 \\
\hline
 6& 406.29& 335.00& 130.00& 130.00& 150.00& 460.00& 465.00&  60.00&  25.00&  25.00&  42.50&  53.13&  25.00&  15.00&  15.00&  20.92 \\
\hline
 7& 416.88& 337.78& 130.00& 130.00& 150.00& 460.00& 465.00&  60.00&  25.00&  25.00&  42.50&  55.00&  25.00&  15.00&  15.00&  21.17 \\
\hline
 8& 455.00& 414.00& 130.00& 130.00& 150.00& 460.00& 465.00&  60.00&  25.00&  25.00&  42.50&  55.00&  25.00&  15.00&  15.00&  23.55 \\
\hline
 9& 455.00& 455.00& 130.00& 130.00& 230.00& 460.00& 465.00&  60.00&  25.00&  53.78&  80.00&  80.00&  25.00&  15.00&  15.00&  27.83 \\
\hline
10& 455.00& 455.00& 130.00& 130.00& 305.00& 460.00& 465.00&  60.00&  25.00&  59.15&  80.00&  80.00&  25.00&  15.00&  15.00&  31.14 \\
\hline
11& 455.00& 455.00& 130.00& 130.00& 346.13& 460.00& 465.00&  60.00&  25.00&  75.63&  80.00&  80.00&  25.00&  15.00&  15.00&  33.76 \\
\hline
12& 455.00& 455.00& 130.00& 130.00& 348.26& 460.00& 465.00&  60.00&  25.00&  75.63&  80.00&  80.00&  25.00&  15.00&  15.00&  33.88 \\
\hline
13& 455.00& 455.00& 130.00& 130.00& 342.94& 460.00& 465.00&  60.00&  25.00&  75.63&  80.00&  80.00&  25.00&  15.00&  15.00&  33.58 \\
\hline
14& 455.00& 455.00& 130.00& 130.00& 390.00& 460.00& 465.00&  60.00&  25.00&  81.57&  80.00&  80.00&  25.00&  15.00&  15.00&  36.57 \\
\hline
15& 455.00& 455.00& 130.00& 130.00& 470.00& 460.00& 465.00&  60.00&  25.00& 132.24&  80.00&  80.00&  25.00&  15.00&  15.00&  44.27 \\
\hline
16& 455.00& 455.00& 130.00& 130.00& 470.00& 460.00& 465.00&  60.00&  25.00& 129.14&  80.00&  80.00&  25.00&  15.00&  15.00&  44.12 \\
\hline
17& 455.00& 455.00& 130.00& 130.00& 438.54& 460.00& 465.00&  60.00&  25.00& 109.38&  80.00&  80.00&  25.00&  15.00&  15.00&  40.89 \\
\hline
18& 455.00& 455.00& 130.00& 130.00& 379.33& 460.00& 465.00&  60.00&  25.00&  63.98&  80.00&  80.00&  25.00&  15.00&  15.00&  35.40 \\
\hline
19& 455.00& 455.00& 130.00& 130.00& 259.33& 460.00& 465.00&  60.00&  25.00&  25.00&  80.00&  80.00&  25.00&  15.00&  15.00&  28.42 \\
\hline
20& 455.00& 455.00& 130.00& 130.00& 200.00& 460.00& 465.00&  60.00&  25.00&  25.00&  72.50&  77.63&  25.00&  15.00&  15.00&  26.15 \\
\hline
21& 455.00& 402.79& 130.00& 130.00& 150.00& 460.00& 465.00&  60.00&  25.00&  25.00&  42.50&  55.00&  25.00&  15.00&  15.00&  23.29 \\
\hline
22& 402.22& 335.00& 130.00& 130.00& 150.00& 460.00& 465.00&  60.00&  25.00&  25.00&  42.50&  53.13&  25.00&  15.00&  15.00&  20.85 \\
\hline
23& 390.26& 295.00& 130.00& 130.00& 150.00& 460.00& 465.00&  60.00&  25.00&  25.00&  42.50&  53.13&  25.00&  15.00&  15.00&  19.90 \\
\hline
24& 383.14& 295.00& 130.00& 130.00& 150.00& 460.00& 465.00&  60.00&  25.00&  25.00&  42.50&  53.13&  25.00&  15.00&  15.00&  19.78 \\
\hline
\end{tabular}}
\end{table}

\section{Conclusion}
\label{Conclusion}
In this paper, an MILP method has been successfully introduced to solve the DED problem considering POZ, ramp rate constraints, transmission losses and spinning reserve constraints.
By using the perspective cut reformulation and first order Taylor expansion techniques, two MILP formulations, i.e., MILP-1 and MILP-2, are formed for the two cases of the DED-POZ.
Based on these two MILP formulations, a straightforward MILP-IA is proposed to optimize the DED-POZ. In order to assess the quality of the solutions, not only the optimality but also the feasibility are discussed. The simulation results show that both MILP-1 and MILP-IA can solve to competitive solutions in a short time. In other words, the proposed
MILP method can solve the DED-POZ problem efficiently.

\section*{Acknowledgment}

This work was supported by the Natural Science Foundation of China (No.514\\
07037); the Natural Science Foundation of Guangxi (No.2014GXNSFFA118001); the Open Project Program of Guangxi Colleges and Universities Key Laboratory of Complex System Optimization and Big Data Processing (No.2016CSOBDP02\\
05).

\end{document}